\documentclass{article}

\usepackage{amsthm}
\usepackage{amssymb}
\usepackage{amsxtra}

\newtheorem{thm}{Theorem}[section]
\newtheorem{lemma}{Lemma}[section]

\title{Lower order terms of the second moment of $S(t)$}
\author{Tsz Ho Chan}

\begin{document}
\maketitle
\begin{abstract}
In this article, we derive an asymptotic formula for the second moment of $S(t)$ which includes the lower order terms using a prediction from the ratios conjecture of Conrey, Farmer and Zirnbauer. The formula matches very well with actual values of the second moment.
\end{abstract}
\section{Introduction}
Let $\rho = \beta + i\gamma$ denote the non-trivial zeros of the
Riemann zeta function $\zeta(s)$. For $t \neq \gamma$,
\begin{equation}
\label{1}
S(t) := \frac{1}{\pi} \mbox{arg} \zeta(\frac{1}{2} + it) = 
\frac{1}{\pi} \mbox{Im} \log{\zeta(\frac{1}{2} + it)},
\end{equation}
where the argument is obtained by continuous variation along the 
horizontal line $\sigma + it$ starting with value zero at $\infty
+ it$. For $t = \gamma$, we define
$$S(t) := \lim_{\epsilon \rightarrow 0} \frac{S(t + \epsilon) - S(t -
\epsilon)}{2}.$$
Selberg [\ref{S}] proved that
$$\int_{0}^{T} |S(t)|^2 dt = \frac{T}{2\pi^2} \log \log T + O(T (\log
\log T)^{1/2}).$$
He also proved the same asymptotic formula with a better error term
$O(T)$ under Riemann Hypothesis ({\bf RH}). Later, Goldston [\ref{G}]
proved that, under {\bf RH},
\begin{equation}
\label{Goldston}
\begin{split}
\int_{0}^{T} |S(t)|^2 dt =& \frac{T}{2\pi^2} \log{\log{\frac{T}{2\pi}}}
+ \frac{T} {2\pi^2} \Bigl[ \int_{1}^{\infty} \frac{F(\alpha, T)}{\alpha^2}
d\alpha + C_0 \\
&- \sum_{m=2}^{\infty} \sum_p \Bigl(\frac{1}{m} - \frac{1}{m^2}
\Bigr) \frac{1}{p^m} \Bigr] + o(T)
\end{split}
\end{equation}
where $C_0$ is Euler's constant and
$$F(\alpha, T) = \Bigl(\frac{T}{2\pi} \log{\frac{T}{2\pi}}\Bigr)^{-1}
\sum_{0 < \gamma, \gamma' \leq T} \Bigl(\frac{T}{2\pi}\Bigr)^{i \alpha
(\gamma - \gamma')} w(\gamma - \gamma'),$$
with $w(u) = 4/(4+u^2)$, is Montgomery's pair correlation function. Here
and throughout the paper, $p$ will denote a prime and sums or products 
over $p$ are over all primes. Note that the integral involving $F$ equals
to $1 + o(1)$ if Montgomery's Strong pair correlation conjecture is true.
Under {\bf RH} and a strong quantitative form of Twin Prime Conjecture,
the author proved in [\ref{C}] that
\begin{equation*}
\begin{split}
\int_{0}^{T} |S(t)|^2 dt =& \frac{T}{2\pi^2} \log{\log{\frac{T}{2\pi}}}
+ \frac{T} {2\pi^2} \Bigl[ \int_{1}^{\infty} \frac{F(\alpha, T)}{\alpha^2}
d\alpha + C_0 \\
&- \sum_{m=2}^{\infty} \sum_p \Bigl(\frac{1}{m} - \frac{1}{m^2}
\Bigr) \frac{1}{p^m} \Bigr] + \frac{C_1 T}{2\pi^2 \log^2{(T/2\pi)}} +
O\Bigl(\frac{T \log \log T}{\log^3 T}\Bigr)
\end{split}
\end{equation*}
where $C_1$ is some constant. This does not quite give the next lower order term $\frac{T}{\log^2 (T/2\pi)}$ as things are hidden in the integral involving $F$.
Our goal is to get all the lower order terms for the second moment. To
this, we shall prove
\begin{thm}
\label{theorem}
Assume {\bf RH} and (\ref{ratio}), we have
\begin{equation}
\label{main}
\begin{split}
\int_{0}^{T} & |S(t)|^2 dt = \frac{1}{2\pi^2} \int_{0}^{T} \log{\log{\frac{t}{2\pi}}}
+ \Bigl[1 + C_0 - \sum_{m=2}^{\infty} \sum_p \Bigl(\frac{1}{m} - \frac{1}{m^2}\Bigr) \frac{1}{p^m} \Bigr] dt\\
&+ \frac{1}{\pi} \int_{0}^{1} \Bigl[ \zeta(1 - v) \zeta(1 + v) A(v) + \frac{1}{v^2} \Bigr] \frac{(\frac{T}{2\pi})^{1-v} - x^{1-v}}{1-v} v \; dv \\
&+ O(T e^{-c (\log x)^{\alpha}})
\end{split}
\end{equation}
where $A$ is defined in (\ref{a}), $0 < \alpha < 1$, $0 < c < 1$ may depend on $\alpha$, and $x = T^{\beta}$ with free parameter $0 < \beta < 1/2$.
\end{thm}
Although it is not immediately clear why the second integral above gives the lower order terms, equation (\ref{II}) in the proof gives some indications. The last section shows very good numerical evidence for the above Theorem!

One can modify the proof to get short interval result (left for the readers):
\begin{thm}
Assume {\bf RH} and (\ref{ratio1}), we have
\begin{equation*}
\begin{split}
\int_{T}^{T+H} & |S(t)|^2 dt = \frac{1}{2\pi^2} \int_{T}^{T+H} \log{\log{\frac{t}{2\pi}}}
+ \Bigl[1 + C_0 - \sum_{m=2}^{\infty} \sum_p \Bigl(\frac{1}{m} - \frac{1}{m^2}\Bigr) \frac{1}{p^m} \Bigr] dt\\
&+ \frac{1}{\pi} \int_{0}^{1} \Bigl[ \zeta(1 - v) \zeta(1 + v) A(v) + \frac{1}{v^2} \Bigr] \frac{(\frac{T+H}{2\pi})^{1-v} - (\frac{T}{2\pi})^{1-v}}{1-v} v \; dv \\
&+ O(H e^{-c (\log T)^{\alpha}})
\end{split}
\end{equation*}
for $T^\epsilon \leq H \leq T$ with $\epsilon > 0$. Again $0 < \alpha < 1$ and $0 < c < 1$ may depend on $\alpha$.
\end{thm}

{\bf Acknowledgement} The author would like to thank Brian Conrey for suggesting the problem and many helpful discussions. Also, he thanks Matthew Young, Sidney Graham, Michael Rubinstein and Daniel Goldston for valuable suggestions. Lastly, the author thanks the American Institute of Mathematics for providing a stimulating environment to work at.
\section{Explicit formulas}
We need some notations for inverse Mellin transform. For real $f(x)$, its
Mellin transform is
$$F(s) := \int_{0}^{\infty} x^s f(x) \frac{dx}{x}.$$
We define
$$\tilde{F}(x) := \frac{1}{2\pi i} \int_{c-i \infty}^{c+i \infty}  F(s)
x^{-s} ds$$
for some $c$ such that the integral converges. Here $\tilde{F}$ is called
the inverse Mellin transform of $F$ and $\tilde{F}(x) = f(x)$ if $f(x)$ is
a nice function.

We need a certain nice entire function $g(z)$. From Gelfand and Shilov
[\ref{GS}], Chap IV, Sect. 8, there is an entire function $f(z) \not\equiv 0$
satisfying
$$|f(x + iy)| \ll e^{|y|/4 - |x|^{\alpha}}$$
for $0 < \alpha < 1$, and $f(z)$ is real for real $z$ (i.e. $\overline{f(z)} = f(\overline{z})$). Without loss of generality, we can assume $f(0) = 1$.
Define $f_1(z) = \frac{f(z) + f(-z)}{2}$. Then $f_1(z) = f_1(-z)$, $\overline{f_1(z)} = f_1(\overline{z})$, and
$$|f_1(x + iy)| \ll e^{|y|/4 - |x|^{\alpha}}.$$
Define $f_2(z) = f_1(iz)$. We have
$$|f_2(x + iy)| \ll e^{|x|/4 - |y|^{\alpha}}$$
and, for real $\sigma$,
\begin{equation}
\label{con}
\overline{f_2(\sigma)} = \overline{f_1(i\sigma)} = f_1(-i\sigma) =
f_1(i\sigma) = f_2(\sigma).
\end{equation}
Finally, we define an entire function $g(z) := e^{z/2} f_2(z)$. By (\ref{con}),
$g(z)$ is real when $z$ is real (i.e. $\overline{g(z)} = g(\overline{z})$).
Moreover,
\begin{equation}
\label{ineq}
g(x+iy) \ll \Bigl\{ \begin{array}{ll}
e^{x - |y|^{\alpha}}, & \mbox{ if } x \geq 0; \\
e^{x/4 - |y|^{\alpha}}, & \mbox{ if } x < 0.
\end{array}
\end{equation}
In summary, $g(z)$ is an entire function of order $1$ with good decay in the negative $x$ direction and $y$-directions, and $g(0) = 1$.
\begin{lemma}
\label{lemma1}
For $x \geq 4$ and $s \neq 1, \; \rho, \; -2n$,
\begin{equation*}
\begin{split}
\frac{\zeta'}{\zeta}(s) =& - \sum_{n=1}^{\infty} \frac{\Lambda(n)}{n^s} \tilde{g_0} (e^{\frac{\log n}{\log x}}) + \frac{g((1-s)\log x)}{1-s} - \sum_{\rho}
\frac{g((\rho - s) \log x)}{\rho -s} \\
&+ \sum_{n=1}^{\infty} \frac{g(-(2n+s)\log x)}{2n + s}
\end{split}
\end{equation*}
where $g_0(z) = g(z)/z$ and $g(z)$ is the entire function constructed above.
\end{lemma}

Proof: Consider the integral
$$I = \frac{1}{2\pi i} \int_{c-i \infty}^{c+i \infty} \frac{\zeta'}{\zeta}(s+w) \frac{g(w \log x)}{w} dw$$
where $c > 1 - \Re s$. On one hand, by Cauchy's residue theorem after shifting the line of integration to the left allowed by (\ref{ineq}), and $g(0) = 1$, we have
\begin{equation}
\label{1.1}
I = \frac{\zeta'}{\zeta}(s) - \frac{g((1-s)\log x)}{1-s} + \sum_{\rho} \frac{g((\rho - s)\log x)}{\rho - s} + \sum_{n = 1}^{\infty} \frac{g(-(2n - s) \log x)}{-(2n + s)}
\end{equation}
due to the pole of $\frac{1}{w}$ at $w=0$, the pole of $\frac{\zeta'}{\zeta}(s+w)$ at $w = 1-s$, and the poles of the zeros of $\zeta(s+w)$ at $w = \rho-s$ and $w = -2n-s$ respectively. Note that the sums over $\rho$ and $n$ converge because of the nice decay (\ref{ineq}) of $g$. On the other hand, using $\frac{\zeta'}{\zeta}(s) = - \sum_n \frac{\Lambda(n)}{n^s}$,
\begin{equation}
\label{1.2}
\begin{split}
I =& - \sum_{n=1}^{\infty} \frac{\Lambda(n)}{n^s} \frac{1}{2\pi i} \int_{c-i \infty}^{c+i \infty} \frac{g(w \log x)}{w} n^{-w} dw \\
=& - \sum_{n=1}^{\infty} \frac{\Lambda(n)}{n^s} \frac{1}{2\pi i} \int_{c-i \infty}^{c+i \infty} \frac{g(v)}{v} e^{-\frac{\log n}{\log x} v}dv
= - \sum_{n=1}^{\infty} \frac{\Lambda(n)}{n^s} \tilde{g_0}(e^{\frac{\log n}
{\log x}}).
\end{split}
\end{equation}
Here the interchange of summation and integration and the convergence of the last sum are justified by the good decay of $g$ even when $\Re s \leq 1$. Combining (\ref{1.1}) and (\ref{1.2}), we have the lemma.
\begin{lemma}
\label{lemma2}
Assume {\bf RH}. For $t \geq 1$, $t \neq \gamma$, $x \geq 4$, we have
\begin{equation*}
\begin{split}
S(t) =& -\frac{1}{\pi} \sum_{n=1}^{\infty} \frac{\Lambda(n) \sin{(t \log n)}}{n^{1/2} \log n} \tilde{g_0}(e^{\frac{\log n}{\log x}}) \\
&+ \sum_{\gamma} h((t - \gamma)\log x) + O\Bigl(\frac{x^{1/2}}{t} e^{-|t \log x|^{\alpha}} \Bigr) + O\Bigl(\frac{e^{-|t \log x|^{\alpha}}}{x^{1/2}} \Bigr)
\end{split}
\end{equation*}
where
\begin{equation}
\label{h}
h(v) := \frac{1}{2\pi i} \int_{-\infty}^{0} \frac{g(u - iv)}{u - iv} - \frac{g(u + iv)}{u + iv} du.
\end{equation}
\end{lemma}

Proof: By Lemma \ref{lemma1} with $s = \sigma + it$, $\sigma > 0$,
\begin{equation*}
\begin{split}
\mbox{Im} \frac{\zeta'}{\zeta}(\sigma + it) =& \sum_{n=1}^{\infty} \frac{\Lambda(n)}{n^\sigma} \sin{(t \log n)} \tilde{g_0} (e^{\frac{\log n}{\log x}}) + \mbox{Im} \frac{g((1-\sigma - it)\log x)}{1-\sigma - it} \\
&- \mbox{Im} \sum_{\rho} \frac{g((\rho - \sigma - it) \log x)}{\rho - \sigma - it} + \mbox{Im} \sum_{n=1}^{\infty} \frac{g(-(2n + \sigma + it)\log x)}{2n + \sigma + it}. \end{split}
\end{equation*}
By (\ref{ineq}),
$$\sum_{n=1}^{\infty} \frac{g(-(2n + \sigma + it)\log x)}{2n + \sigma + it} \ll \sum_{n=1}^{\infty} e^{-(n/2 + \sigma/4)\log x} e^{-|t \log x|^{\alpha}} \ll \frac{e^{-|t \log x|^{\alpha}}}{x^{1/2 + \sigma/4}}.$$
Then, by (\ref{1}),
\begin{equation}
\label{1.3}
\begin{split}
& S(t) = - \frac{1}{\pi} \int_{1/2}^{\infty} \mbox{Im} \frac{\zeta'}{\zeta}(\sigma + it) d\sigma \\
=& - \frac{1}{\pi} \sum_{n=1}^{\infty} \frac{\Lambda(n) \sin{(t \log n)}}{n^{1/2} \log n} \tilde{g_0} (e^{\frac{\log n}{\log x}}) + \frac{1}{\pi} \int_{1/2}^{\infty} \mbox{Im} \sum_{\rho} \frac{g((\rho - \sigma - it) \log x)}{\rho - \sigma - it} d\sigma \\
&+ O\Bigl(\frac{x^{1/2}}{t} e^{-|t \log x|^{\alpha}} \Bigr) + O\Bigl(\frac{e^{-|t \log x|^{\alpha}}}{x^{1/2}} \Bigr)
\end{split}
\end{equation}
because, by (\ref{ineq}),
\begin{equation*}
\begin{split}
& \int_{1/2}^{\infty} \frac{g((1 - \sigma - it)\log x)}{1 - \sigma -it} d\sigma \\
\ll& \frac{1}{t} \Bigl[\int_{1/2}^{1} e^{(1-\sigma)\log x} e^{-|t \log x|^{\alpha}} d\sigma + \int_{1}^{\infty} e^{-\frac{(\sigma-1)\log x}{4}} e^{-|t \log x|^{\alpha}} d\sigma \Bigr]\\
\ll& \frac{x^{1/2}}{t} e^{-|t \log x|^{\alpha}}.
\end{split}
\end{equation*}
Now, using {\bf RH}, $\mbox{Im} z = \frac{1}{2i}(z - \overline{z})$ and $\overline{g(z)} = g(\overline{z})$, the integral in (\ref{1.3})
\begin{equation}
\label{1.4}
\begin{split}
=& \int_{1/2}^{\infty} \mbox{Im} \sum_{\gamma} \frac{g(((\frac{1}{2} - \sigma) - i(t - \gamma)) \log x)}{(\frac{1}{2} - \sigma) - i(t - \gamma)} d\sigma \\
=& \sum_{\gamma} \frac{1}{2i} \int_{-\infty}^{0} \frac{g((\sigma - i(t - \gamma)) \log x)}{\sigma - i(t - \gamma)} - \frac{g((\sigma + i(t - \gamma)) \log x)}{\sigma + i(t - \gamma)} d\sigma \\
=& \sum_{\gamma} \frac{1}{2i} \int_{-\infty}^{0} \frac{g(u - i(t - \gamma) \log x)}{u - i(t - \gamma) \log x} - \frac{g(u + i(t - \gamma) \log x)}{u + i(t - \gamma)\log x} du.
\end{split}
\end{equation}
Again, the interchange of summation and integration and the convergence of the final sum are justified by the good decay of $g$. Hence, the lemma follows from (\ref{1.3}) and (\ref{1.4}).
\section{Some properties of $h(v)$}
Recall $h$ is given by (\ref{h}). Clearly, $h(0) = 0$ and $h(v)$ is an odd, real-valued function. For the decay of $h$, we have
\begin{lemma}
\label{lemma3}
$$h(v) \ll e^{-|v|^{\alpha}}.$$
\end{lemma}

Proof: Since $g(z)$ is entire and $g(0) = 1$, $g(z) = 1 + O(z)$ when $|z| \leq \epsilon$ for some small fixed $\epsilon > 0$. Now, if $|v| > \epsilon$, we have
$$h(v) = \frac{1}{2\pi i} \int_{-\infty}^{0} \frac{g(u-iv)}{u-iv} - \frac{g(u+iv)}{u+iv} du \ll \frac{1}{\epsilon} \int_{-\infty}^{0} e^{-u/4} e^{-|v|^{\alpha}} du \ll e^{-|v|^{\alpha}}.$$
If $0 < |v| \leq \epsilon$, it suffice to consider $0 < v \leq \epsilon$ as $h$ is odd. Observe that $h(v)$ consists of two horizontal line integrals of $\frac{g(z)}{z}$, one at height $-iv$ and the other at height $iv$. So, if we join the two lines by a semi-circle $C_v:$ $\{z = v e^{i\theta}, -\pi/2 \leq \theta \leq \pi/2\}$, we can apply Cauchy's residue theorem and get
\begin{equation*}
\begin{split}
1 =& g(0) = h(v) + \frac{1}{2\pi i} \int_{C_v} \frac{g(z)}{z} dz \\
=& h(v) + \frac{1}{2\pi i} \int_{-\pi/2}^{\pi/2} \frac{1 + O(v)}{v e^{i\theta}} v i e^{i \theta} d\theta \\
=& h(v) + \frac{1}{2\pi} \int_{-\pi/2}^{\pi/2} 1 + O(v) d\theta.
\end{split}
\end{equation*}
Hence $h(v) = \frac{1}{2} + O(v)$ and we have the lemma.
\begin{lemma}
\label{lemma4}
Let $\hat{f}(a) = \int_{-\infty}^{\infty} f(v) e^{-2\pi i a v} dv$ be the Fourier transform of $f$. Then, for large $M, N > 0$ and small enough $\epsilon > 0$,
\begin{equation*}
\begin{split}
\hat{h}(a) =& \frac{1}{2\pi i} \int_{-N}^{-\epsilon^2} \int_{-M}^{M} \Bigl[\frac{g(u - iv)}{u - iv} - \frac{g(u + iv)}{u + iv} \Bigr] e^{-2\pi i a v} dv \; du \\
&+ O(Me^{-M^{\alpha}}) + O(e^{-N/4}) + O(\epsilon).
\end{split}
\end{equation*}
\end{lemma}

Proof: By the definition of $h$ and Lemma \ref{lemma3}, $\hat{h}(a)$ is clearly well-defined and
\begin{equation}
\label{2.1}
\hat{h}(a) = \int_{-M}^{M} h(v) e^{-2\pi i a v} dv + O(M e^{-M^{\alpha}}).
\end{equation}
Putting in the definition of $h(v)$ and using (\ref{ineq}), we have
\begin{equation}
\label{2.2}
\begin{split}
&\int_{-M}^{M} h(v) e^{-2\pi i a v} dv \\
=& \frac{1}{2\pi i} \int_{-M}^{M} \int_{-\infty}^{0} \Bigl[\frac{g(u - iv)}{u - iv} - \frac{g(u + iv)}{u+iv} \Bigr] du \; e^{-2\pi i a v} dv \\
=& \frac{1}{2\pi i} \int_{-M}^{M} \int_{-N}^{0} \Bigl[\frac{g(u - iv)}{u - iv} - \frac{g(u + iv)}{u+iv} \Bigr] du \; e^{-2\pi i a v} dv  + O(e^{-N/4}) \\
=& \frac{1}{2\pi i} \int_{-M}^{M} \int_{-N}^{-\epsilon^2} \Bigl[\frac{g(u - iv)}{u - iv} - \frac{g(u + iv)}{u+iv} \Bigr] du \; e^{-2\pi i a v} dv \\
+& \frac{1}{2\pi i} \int_{-M}^{M} \int_{-\epsilon^2}^{0} \Bigl[\frac{g(u - iv)}{u - iv} - \frac{g(u + iv)}{u+iv} \Bigr] du \; e^{-2\pi i a v} dv + O(e^{-N/4}) \\
=& I_1 + I_2 + O(e^{-N/4}).
\end{split}
\end{equation}
$I_2$ can be broken down further into
$$\frac{1}{2\pi i} \Bigl[ \int_{\epsilon}^{M} \int_{-\epsilon^2}^{0} + \int_{-M}^{-\epsilon} \int_{-\epsilon^2}^{0} + \int_{-\epsilon}^{\epsilon} \int_{-\epsilon^2}^{0} \Bigr] = \frac{1}{2\pi i} [ J_1 + J_2 + J_3 ].$$
By (\ref{ineq}),
\begin{equation*}
\begin{split}
J_1 =& \int_{-\epsilon^2}^{0} \int_{\epsilon}^{M} \Bigl[\frac{g(u - iv)}{u - iv} - \frac{g(u + iv)}{u+iv} \Bigr] e^{-2 \pi i a v} dv \; du \\
\ll& \int_{-\epsilon^2}^{0} \frac{1}{\epsilon} \int_{\epsilon}^{M} e^{-u/4} e^{-|v|^{\alpha}} dv \; du \ll \int_{-\epsilon^2}^{0} \frac{1}{\epsilon} e^{-u/4} du \ll \epsilon.
\end{split}
\end{equation*}
Similarly,
$$J_2 \ll \epsilon.$$
Now, as $g$ is entire and $g(0) = 1$, $g(u \pm iv) = 1 + O(|u \pm iv|)$ for $-\epsilon^2 \leq u \leq 0$ and $-\epsilon \leq v \leq \epsilon$ provided $\epsilon$ is small enough. Thus, by oddness of the inner integral in $v$,
\begin{equation*}
\begin{split}
J_3 =& -i \int_{-\epsilon}^{\epsilon} \int_{-\epsilon^2}^{0} \Bigl[\frac{g(u - iv)}{u - iv} - \frac{g(u + iv)}{u+iv} \Bigr] \sin{(2\pi a v)} du \; dv \\
=& -i \int_{-\epsilon}^{\epsilon} \int_{-\epsilon^2}^{0} \Bigl[ \frac{1}{u-iv} - \frac{1}{u+iv} + O(1) \Bigr] \sin{(2\pi a v)} du \; dv \\
=& -i \int_{-\epsilon}^{\epsilon} \int_{-\epsilon^2}^{0} \frac{2 v \sin{(2\pi a v)}}{u^2 + v^2} du \; dv + O(\epsilon^3) \ll \epsilon^3.
\end{split}
\end{equation*}
as the integrand is bounded. Therefore, $I_2 \ll \epsilon$. Hence, combining this with (\ref{2.1}) and (\ref{2.2}), we have the lemma.
\begin{lemma}
\label{lemma5}
$$\hat{h}(a) = -\frac{i}{2\pi a} \mbox{ if } |a| > \frac{1}{2\pi}.$$
\end{lemma}

Proof: Since $\hat{h}$ is odd, we may assume $a > \frac{1}{2\pi}$. From Lemma \ref{lemma4}, with $E = Me^{-M^{\alpha}} + e^{-N/4} + \epsilon$,
\begin{equation*}
\begin{split}
\hat{h}(a) =& \frac{1}{2\pi i} \int_{-N}^{-\epsilon^2} \int_{-M}^{M} \frac{g(u - iv)}{u - iv} e^{-2\pi i a v} dv \; du \\
&- \frac{1}{2\pi i} \int_{-N}^{-\epsilon^2} \int_{-M}^{M} \frac{g(u + iv)}{u + iv} e^{-2\pi i a v} dv \; du + O(E) \\
=& \int_{-N}^{-\epsilon^2} e^{-2\pi a u} \frac{1}{2\pi i} \int_{-M}^{M} \frac{g(u - iv)}{u - iv} e^{2\pi a(u -i v)} dv \; du \\
&- \int_{-N}^{-\epsilon^2} e^{2\pi a u} \frac{1}{2\pi i} \int_{-M}^{M} \frac{g(u+iv)}{u+iv} e^{-2\pi a(u + i v)} dv \; du + O(E) \\
=& I_1 - I_2 + O(E).
\end{split}
\end{equation*}
The inner integral of $I_1$
$$= -\frac{1}{i} \frac{1}{2\pi i} \int_{u-iM}^{u+iM} \frac{g(z)}{z} e^{2\pi a z} dz = -\frac{1}{i} \Bigl[0 + O(\frac{1}{M})\Bigr]$$
by shifting the line of integration to the left as $a>0$ and using Cauchy's residue theorem. The error term comes from estimating the two horizontal line integrals at heights $-iM$ and $iM$. Similarly, the inner integral of $I_2$
$$= \frac{1}{i} \frac{1}{2\pi i} \int_{u-iM}^{u+iM} \frac{g(z)}{z} e^{-2\pi a z} dz = \frac{1}{i} \Bigl[-1 + O(\frac{1}{M})\Bigr]$$
by shifting the line of integration to the right and picking up a pole at $z=0$. The shifted line integrals are small provided that $2\pi a > 1$ as $g(z) \ll e^{|z|}$ by (\ref{ineq}). Consequently, for $a > \frac{1}{2\pi}$,
$$\hat{h}(a) = \frac{1}{i} \int_{-N}^{-\epsilon^2} e^{2\pi a u} du + O(\frac{1}{M} + E) = -\frac{i}{2\pi a} [e^{-2\pi a \epsilon^2} - e^{-2\pi a N}] + O(\frac{1}{M} + E).$$
The lemma follows by taking $M, N \rightarrow \infty$ and $\epsilon \rightarrow 0$.
Note: This property of $\hat{h}$ is the same as that of equation (3.10) in Goldston [\ref{G}].

Lastly, we have a crucial identity relating $\tilde{g_0}$ and $\hat{h}$.
\begin{lemma}
\label{lemma6}
$$[1 - \tilde{g_0}(e^u)]^2 = - \Bigl(u \hat{h}\Bigl(\frac{u}{2\pi}\Bigr) \Bigr)^2.$$
\end{lemma}

Proof: Recall the definition of $g_0$ in Lemma \ref{lemma1}, it suffices to show
\begin{equation}
\label{2.3}
1 - \frac{1}{2\pi i} \int_{c-i\infty}^{c+i\infty} \frac{g(v)}{v} e^{-uv} dv = i u \hat{h}\Bigl(\frac{u}{2\pi}\Bigr)
\end{equation}
for $c > 0$. Firstly, one can easily check that (\ref{2.3}) is true when $u=0$ by Cauchy's residue theorem and shifting of the line of integration. Without loss of generality, we may assume $u > 0$ as the other case is similar. Now, by Lemma \ref{lemma4}, it suffices to show, with $E = Me^{-M^{\alpha}} + e^{-N/4} + \epsilon$,
\begin{equation}
\label{2.4}
\begin{split}
1 =& \frac{1}{2\pi i} \int_{c-i\infty}^{c+i\infty} \frac{g(z)}{z} e^{-uz} dz + \frac{i u}{2\pi i} \int_{-N}^{-\epsilon^2} \int_{-M}^{M} \frac{g(\sigma -iv)}{\sigma - iv} e^{-i u v} dv \; d\sigma \\
&- \frac{i u}{2\pi i} \int_{-N}^{-\epsilon^2} \int_{-M}^{M} \frac{g(\sigma + iv)}{\sigma + iv} e^{-i u v} dv \; d\sigma + O(uE) \\
=& I_1 + I_2 + I_3 + O(uE).
\end{split}
\end{equation}
Imitating the calculation in Lemma \ref{lemma5}, one has
$$I_2 = -u \int_{-N}^{-\epsilon^2} e^{-\sigma u} \frac{1}{2\pi i} \int_{\sigma - iM}^{\sigma + iM} \frac{g(z)}{z} e^{uz} dz \; d\sigma = O(\frac{1}{M});$$
and
$$I_3 = -u \int_{-N}^{-\epsilon^2} e^{\sigma u} \frac{1}{2\pi i} \int_{\sigma - iM}^{\sigma + iM} \frac{g(z)}{z} e^{-uz} dz \; d\sigma.$$
Instead of shifting of the line of integration to the right in $I_3$, we just shift it to the vertical line $\Re z = -\epsilon$ with an error $O(\frac{1}{M})$. Thus,
\begin{equation*}
\begin{split}
I_3 =& -u \frac{1}{2\pi i} \int_{-\epsilon-iM}^{-\epsilon+iM} \frac{g(z)}{z} e^{-uz} dz \int_{-N}^{-\epsilon^2} e^{\sigma u} d\sigma + O(\frac{1}{M}) \\
=& - \frac{1}{2\pi i} \int_{-\epsilon-iM}^{-\epsilon+iM} \frac{g(z)}{z} e^{-uz} dz (1 + O(\epsilon^2 u + e^{-N u})) + O(\frac{1}{M}) \\
=& - \frac{1}{2\pi i} \int_{-\epsilon-i\infty}^{-\epsilon+i\infty} \frac{g(z)}{z} e^{-uz} dz (1 + O(\epsilon^2 u + e^{-N u})) + O(\frac{1}{M}) \\
=& - \frac{1}{2\pi i} \int_{-\epsilon-i\infty}^{-\epsilon+i\infty} \frac{g(z)}{z} e^{-uz} dz + O\Bigl(\epsilon u + \frac{e^{-N u}}{\epsilon} + \frac{1}{M} \Bigr)
\end{split}
\end{equation*}
by (\ref{ineq}) again. Therefore, the right hand side of (\ref{2.4})
$$= \frac{1}{2\pi i} \int_{c-i\infty}^{c+i\infty} \frac{g(z)}{z} e^{-uz} dz - \frac{1}{2\pi i} \int_{-\epsilon-i\infty}^{-\epsilon+i\infty} \frac{g(z)}{z} e^{-uz} dz + O\Bigl(\epsilon u + \frac{e^{-N u}}{\epsilon} + \frac{1}{M} + uE \Bigr)$$
which equals to $1$ by Cauchy's residue theorem and letting $\epsilon \rightarrow 0$; $M, N \rightarrow \infty$ appropriately. Hence, we have the lemma.

As a consequence of Lemmas \ref{lemma5} and \ref{lemma6}, we have
\begin{equation}
\label{amazing}
\tilde{g_0}(e^u) = 0 \mbox{ when } |u| > 1.
\end{equation}
This can also be proved directly by considering $\frac{1}{2\pi i} \int_{c-i\infty}^{c+i\infty} \frac{g(s)}{s} e^{-u s} ds$ and shifting the line of integration to the right.
\section{Starting the proof of Theorem \ref{theorem}}
The method of proof is essentially that of Goldston [\ref{G}]. By Lemma \ref{lemma2} and (\ref{amazing}), we have, for $t \geq 1$ and $t \neq \gamma$,
\begin{equation*}
\begin{split}
S(t) +& \frac{1}{\pi} \sum_{n \leq x}^{} \frac{\Lambda(n) \sin{(t \log n)}}{n^{1/2} \log n} \tilde{g_0}(e^{\frac{\log n}{\log x}}) \\
=& \sum_{\gamma} h((t - \gamma)\log x) + O\Bigl(\frac{x^{1/2}}{t} e^{-|t \log x|^{\alpha}} \Bigr) + O\Bigl(\frac{e^{-|t \log x|^{\alpha}}}{x^{1/2}} \Bigr).
\end{split}
\end{equation*}
Since $h$ is bounded and the above formula holds except on a countable set of points, we have on squaring both sides and integrating from $1$ to $T$,
\begin{equation}
\label{3.1}
\int_{1}^{T} |S(t)|^2 dt + H(T) + G(T) = R + O(T^{1/2} x^{1/2}),
\end{equation}
where
$$G(T) = \int_{1}^{T} \Big|\frac{1}{\pi} \sum_{n \leq x} \frac{\Lambda(n) \sin{(t \log n)}}{n^{1/2} \log n} \tilde{g_0}(e^{\frac{\log n}{\log x}}) \Big|^2 dt,$$
$$H(T) = \frac{2}{\pi} \int_{1}^{T} S(t) \sum_{n \leq x} \frac{\Lambda(n) \sin{(t \log n)}}{n^{1/2} \log n} \tilde{g_0}(e^{\frac{\log n}{\log x}}) dt,$$
and
$$R = \int_{1}^{T} \Big| \sum_{\gamma} h((t-\gamma) \log x) \Big|^2 dt.$$
The error term is obtained by Cauchy-Schwarz inequality since $R \ll T$. The lower limit of integration may be replaced by zero since $\int_{0}^{1} |S(t)|^2 dt \ll 1$. Following the same calculation as Goldston [\ref{G}], we have, for $x \geq 4$ and $T \geq 2$,
\begin{equation}
\label{R}
R = \frac{1}{\log x} \sum_{0 < \gamma, \gamma' \leq T} \hat{k}((\gamma - \gamma')\log x) + O(\log^3 T)
\end{equation}
where $k(u) = - \hat{h}(u)^2$. Note that since $h$ is odd, $\hat{h}$ is odd. So $k$ and $\hat{k}$ are even functions. Also, from [\ref{G}],
$$G(T) = \frac{T}{2\pi^2} \sum_{n \leq x} \frac{\Lambda^2(n)}{n \log^2 n} \tilde{g_0}^2 (e^{\frac{\log n}{\log x}}) + O(x^2),$$
and
$$H(T) = - \frac{T}{\pi^2} \sum_{n \leq x} \frac{\Lambda^2(n)}{n \log^2 n} \tilde{g_0}(e^{\frac{\log n}{\log x}}) + O(x^{2+\epsilon}).$$
Thus, by Lemma \ref{lemma6},
\begin{equation*}
\begin{split}
& -G(T) - H(T) \\
=& -\frac{T}{2\pi^2} \sum_{n \leq x} \frac{\Lambda^2(n)}{n \log^2 n} \Bigl[ 1 - \tilde{g_0}(e^{\frac{\log n}{\log x}}) \Bigr]^2 + \frac{T}{2\pi^2} \sum_{n \leq x} \frac{\Lambda^2(n)}{n \log^2 n} \\
=& -\frac{T}{2\pi^2} \sum_{n \leq x} \frac{\Lambda^2(n)}{n \log^2 n} \Bigl(\frac{\log n}{\log x}\Bigr)^2 k\Bigl(\frac{\log n}{2\pi \log x}\Bigr) + \frac{T}{2\pi^2} \sum_{n \leq x} \frac{\Lambda^2(n)}{n \log^2 n} \\
=& -\frac{T}{2\pi^2} \frac{1}{\log^2 x} \sum_{n \leq x} \frac{\Lambda^2(n)}{n} k\Bigl(\frac{\log n}{2\pi \log x} \Bigr) + \frac{T}{2\pi^2} \sum_{n \leq x} \frac{\Lambda^2(n)}{n \log^2 n} \\
=& -\frac{T}{2\pi^2} \frac{1}{\log^2 x} \sum_{n \leq x} \frac{\Lambda^2(n)}{n} \int_{-\infty}^{\infty} \log x \; \hat{k}(r \log x) e^{2 \pi i \frac{\log n}{2\pi} r} dr + \frac{T}{2\pi^2} \sum_{n \leq x} \frac{\Lambda^2(n)}{n \log^2 n} \\
=& -\frac{T}{2\pi^2} \frac{1}{\log x} \int_{-\infty}^{\infty} \hat{k}(r \log x) \sum_{n \leq x} \frac{\Lambda^2 (n)}{n^{1+ ir}} dr + \frac{T}{2\pi^2} \sum_{n \leq x} \frac{\Lambda^2(n)}{n \log^2 n}
\end{split}
\end{equation*}
because $\hat{k}$ is even. Now
\begin{equation}
\label{3.2}
\sum_{p \leq u} \frac{1}{p} = \log \log u + C_0 + \sum_{p} \Bigl[\log \Bigl(1 - \frac{1}{p}\Bigr) + \frac{1}{p} \Bigr] + r(u),
\end{equation}
where $r(u) \ll u^{-1/2 + \epsilon}$ under {\bf RH} and the sum over primes in (\ref{3.2}) is equal to $-\sum_{m=2}^{\infty} \sum_p \frac{1}{m p^m}$. Therefore,
\begin{equation}
\label{GH}
\begin{split}
-G(T) - H(T) =& -\frac{T}{2\pi^2} \frac{1}{\log x} \int_{-\infty}^{\infty} \hat{k}(r \log x) \sum_{n \leq x} \frac{\Lambda^2 (n)}{n^{1+ ir}} dr + \frac{T}{2\pi^2} \Bigl[ \log \log x + C_0 \\
&- \sum_{m=2}^{\infty} \sum_p \Bigl(\frac{1}{m} - \frac{1}{m^2}\Bigr) \frac{1}{p^m} \Bigr] + O\Bigl(\frac{T}{x^{1/2 - \epsilon}}\Bigr) + O(x^{2+\epsilon}).
\end{split}
\end{equation}
\section{Using the ratios conjecture}
Denote
$$l = \log \frac{t}{2\pi} \; \mbox{ and } \; \lambda = \log x .$$
In view of (\ref{3.1}) and (\ref{GH}), it remains to deal with $R$
or the double sum
$$\Sigma = \sum_{0 < \gamma, \gamma' \leq T} \hat{k}((\gamma - \gamma')\log x) = \sum_{0 < \gamma, \gamma' \leq T} f(\gamma - \gamma')$$
with $f(r) = \hat{k}(r \log x)$. Note: Since $h(v) \ll e^{-|v|^{\alpha}}$ by Lemma \ref{lemma5} and $\hat{k}(v) = h*h(v)$ (the convolution of $h$ with itself), $\hat{k}(v) \ll e^{-|v|^{\alpha}}$ and $f(r) \ll e^{-|r \log x|^{\alpha}}$. It is here that we make use of the ratios conjecture by Conrey, Farmer and Zirnbauer (see Conjecture 1 of [\ref{CFZ}]). Using it, Conrey and Snaith [\ref{CS}] derived a formula for pair correlation of zeros:
For continuous real, even function $f(r)$ with $f(r) \ll 1/(1+r^2)$,
\begin{equation}
\label{ratio}
\begin{split}
\sum_{0 < \gamma, \gamma' \leq T} f(&\gamma - \gamma') = \frac{1}{(2\pi)^2} \int_{0}^{T} \Bigl(2 \pi f(0) \log{\frac{t}{2\pi}} + \int_{-T}^{T} f(r)
\Bigl( \log^2{\frac{t}{2\pi}} + 2 \Bigl( \Bigl(\frac{\zeta'}{\zeta} \Bigr)'
(1+ir) \\
+& \Bigl(\frac{t}{2\pi}\Bigr)^{-ir} \zeta(1-ir) \zeta(1+ir) A(ir) - B(ir)
\Bigr) \Bigr) dr \Bigr) dt + O(T^{1/2 + \epsilon}); 
\end{split}
\end{equation}
here the integral is to be regarded as a principal value near $r=0$,
\begin{equation}
\label{a}
A(\eta) = \prod_p \frac{(1 - \frac{1}{p^{1+\eta}}) (1 - \frac{2}{p} +
\frac{1}{p^{1+\eta}})}{(1 - \frac{1}{p})^2},
\end{equation}
and
\begin{equation}
\label{b}
B(\eta) = \sum_p \Bigl(\frac{\log p}{p^{1+\eta} - 1}\Bigr)^2.
\end{equation}
Remarks: 1. In [\ref{CS}], $f$ is required to be analytic in a horizontal strip around the real axis but our $f$ is not even differentiable. 2. Also, in [\ref{CS}], $f$ has no dependency on $T$ while ours depends on $x$ which in turn depends on $T$.
3. Equation (\ref{ratio}) was first derived by Bogomolny and Keating [\ref{BoK}] using random matrix method (see [\ref{BeK}] for more details).

One may even conjecture, for $\epsilon > 0$ and $T^\epsilon \leq H \leq T$,
\begin{equation}
\label{ratio1}
\begin{split}
\sum_{T < \gamma, \gamma' \leq T+H} f(&\gamma - \gamma') = \frac{1}{(2\pi)^2} \int_{T}^{T+H} \Bigl(2 \pi f(0) \log{\frac{t}{2\pi}} + \int_{-T}^{T} f(r)
\Bigl( \log^2{\frac{t}{2\pi}} + 2 \Bigl( \Bigl(\frac{\zeta'}{\zeta} \Bigr)'
(1+ir) \\
+& \Bigl(\frac{t}{2\pi}\Bigr)^{-ir} \zeta(1-ir) \zeta(1+ir) A(ir) - B(ir)
\Bigr) \Bigr) dr \Bigr) dt + O(H^{1/2} T^\epsilon).
\end{split}
\end{equation}
By (\ref{ratio}),
\begin{equation*}
\begin{split}
\Sigma =& \frac{1}{(2\pi)^2} \int_{0}^{T} \biggl\{2 \pi f(0) l + \int_{-T}^{T} f(r)
\Bigl[ l^2 - \frac{2}{r^2} + 2\frac{e^{-i r l}}{r^2} \Bigr] dr \\
&+ \int_{-T}^{T} f(r) 2 \Bigl[ \Bigl(\frac{\zeta'}{\zeta} \Bigr)'
(1+ir) + \frac{1}{r^2} - B(ir) \Bigr] dr \\
&+ \int_{-T}^{T} f(r) 2 \Bigl(\frac{t}{2\pi}\Bigr)^{-ir} \Bigl[\zeta(1-ir) \zeta(1+ir) A(ir) - \frac{1}{r^2} \Bigr] dr \biggr\} dt + O(T^{1/2 + \epsilon}) \\
=& \frac{1}{(2\pi)^2} \int_{x}^{T} 2\pi l \Bigl[ f(0) + \int_{-\infty}^{\infty} f(\frac{2\pi}{l}u) \Bigl(1 - \Bigl(\frac{\sin{\pi u}}{\pi u}\Bigr)^2 \Bigr) du \Bigr] dt \\
&+ \frac{T}{2\pi^2} \int_{-\infty}^{\infty} f(r) \Bigl[ \Bigl(\frac{\zeta'}{\zeta} \Bigr)'(1+ir) + \frac{1}{r^2} -B(ir) \Bigr] dr \\
&+ \frac{1}{2\pi^2} \int_{x}^{T} \int_{-\infty}^{\infty} f(r) \Bigl[\zeta(1-ir) \zeta(1+ir) A(ir) - \frac{1}{r^2} \Bigr] e^{-ir l } dr \; dt \\
&+ O(x \log x) + O(T^{1/2 + \epsilon}) \\
=& \Sigma_1 + \Sigma_2 + \Sigma_3 + O(x \log x + T^{1/2 + \epsilon}).
\end{split}
\end{equation*}
The integrals can be extended to $\pm \infty$ with small error because of the good decay of $f$. Also, $\Sigma_1$ is essentially done in Conrey and Snaith [\ref{CS}] using appropriate algebra and substitution. Note that the above matches exactly with the pair correlation function $R(x)$ in [\ref{BeK}]. $\Sigma_1$, $\Sigma_2$ and $\Sigma_3$ correspond to $R_{GUE}(x)$, $R_{c}^1(x)$ and $R_{c}^2(x)$ in [\ref{BeK}] respectively.

Since $h$ is odd, one can easily check that $k(0) = \int_{-\infty}^{\infty} h * h(x) dx = 0$. Now, by Parseval's identity $\int_{-\infty}^{\infty} \hat{f}(x) g(x) dx = \int_{-\infty}^{\infty} f(x) \hat{g}(x) dx$,
\begin{equation*}
\begin{split}
\Sigma_1 =& \frac{1}{(2\pi)^2} \int_{x}^{T} 2\pi l \Bigl[ f(0) + \int_{-\infty}^{\infty} f(\frac{2\pi}{l}u) du - \int_{-\infty}^{\infty} f(\frac{2\pi}{l}u) \Bigl(\frac{\sin{\pi u}}{\pi u}\Bigr)^2 du \Bigr] dt \\
=& \frac{1}{(2\pi)^2} \int_{x}^{T} 2\pi l \Bigl[ \hat{k}(0) + 0 - \int_{-\infty}^{\infty} \frac{l}{2\pi \lambda} k \Bigl(\frac{l v}{2\pi \lambda}\Bigr) \max\{1-|v|,0\} dv \Bigr] dt \\
=& \frac{1}{(2\pi)^2} \int_{x}^{T} 2\pi l \Bigl[ \int_{-\infty}^{\infty} \frac{l}{2\pi \lambda} k \Bigl(\frac{l v}{2\pi \lambda}\Bigr) dv - \int_{-1}^{1} \frac{l}{2\pi \lambda} k \Bigl(\frac{l v}{2\pi \lambda}\Bigr) (1 - |v|) dv \Bigr] dt \\
=& \frac{1}{(2\pi)^2} \int_{x}^{T} 4\pi l \Bigl[ \int_{1}^{\infty} \frac{l}{2\pi \lambda} k \Bigl(\frac{l v}{2\pi \lambda}\Bigr) dv + \int_{0}^{1} \frac{l v}{2\pi \lambda} k \Bigl(\frac{l v}{2\pi \lambda}\Bigr) dv \Bigr] dt
\end{split}
\end{equation*}
As $l \geq \lambda$ and $v \geq 1$, $\frac{lv}{2\pi \lambda} \geq \frac{1}{2\pi}$. Thus, since $k(u) = -\hat{h}(u)^2$, we may apply Lemma \ref{lemma5} and get
\begin{equation}
\begin{split}
\Sigma_1 =& \frac{1}{(2\pi)^2} \int_{x}^{T} 4\pi l \Bigl[ \int_{1}^{\infty} \frac{1}{4 \pi^2 (\frac{l v}{2\pi \lambda})^2} \frac{l}{2\pi \lambda} dv + \int_{\lambda/l}^{1} \frac{1}{4 \pi^2 (\frac{l v}{2\pi \lambda})^2} \frac{l v}{2\pi \lambda} dv \\
&+ \int_{0}^{\lambda/l} \frac{l v}{2\pi \lambda} k \Bigl(\frac{l v}{2\pi \lambda}\Bigr) dv \Bigr] dt \\
=& \frac{1}{(2\pi)^2} \int_{x}^{T} 4\pi l \Bigl[ \frac{2\pi \lambda}{4 \pi^2 l} + \frac{2\pi \lambda}{4\pi^2 l} \log \frac{l}{\lambda} + \frac{2\lambda}{l} \int_{0}^{1/2\pi} k(u) u\; du \Bigr] dt \\
=& \int_{x}^{T} \Bigl[ \frac{\lambda}{2\pi^2} + \frac{\lambda}{2\pi^2} \log \frac{l}{\lambda} + 2\lambda \int_{0}^{1/2\pi} k(u) u\; du \Bigr] dt \\
=& \frac{T \lambda}{2\pi^2} + \frac{T \lambda}{2\pi^2} \log \log \frac{T}{2\pi} - \frac{T\lambda}{2\pi^2} \log \log x - \frac{\lambda}{\pi} \mbox{Li}(\frac{T}{2\pi}) \label{3.2.5} \\
&+ 2T \lambda \int_{0}^{1/2\pi} k(u) u \; du + O(x \lambda).
\end{split}
\end{equation}
Here $\mbox{Li}(x) = \int_{0}^{x} \frac{du}{\log u}$.

\bigskip

Next, it turns out that there are some cancellations among $\Sigma_2$, $G(T)$ and $H(T)$. In view of (\ref{3.1}), (\ref{3.2}) and the decomposition of $\Sigma$, let us consider
\begin{equation}
\label{3.3}
\begin{split}
J =& \frac{1}{\lambda} \Sigma_2 - G(T) - H(T) \\
=& \frac{T}{2\pi^2} \Bigl[ \log \log x + C_0 - \sum_{m=2}^{\infty} \sum_p \Bigl(\frac{1}{m} - \frac{1}{m^2}\Bigr) \frac{1}{p^m} \Bigr] \\
&+ \frac{T}{2\pi^2} \frac{1}{\log x} \int_{-\infty}^{\infty} \hat{k}(r \log x) \Bigl[ \Bigl(\frac{\zeta'}{\zeta}\Bigr)' (1+ir) + \frac{1}{r^2} - \sum_p \Bigl(\frac{\log p}{p^{1+ir} - 1}\Bigr)^2 \\
&- \sum_{n \leq x} \frac{\Lambda^2(n)}{n^{1 + ir}} \Bigr] dr + O\Bigl(\frac{T}{x^{1/2 - \epsilon}}\Bigr) = J_1 + J_2 + O\Bigl(\frac{T}{x^{1/2 - \epsilon}}\Bigr) + O(x^{2+\epsilon})
\end{split}
\end{equation}
by (\ref{GH}) and the definitions of $A$ and $B$ (see (\ref{a}) and (\ref{b})).
Before proceeding, we need some lemmas.
\begin{lemma}
\label{lemma7}
For $x > 4$,
$$\sum_p \Bigl(\frac{\log p}{p^{1+ir} - 1}\Bigr)^2 = \sum_{n \leq x} \frac{\Lambda(n) \log n - \Lambda^2(n)}{n^{1+ir}} + O\Bigl(\frac{\log^2 x \log \log x}{x^{1/2}}\Bigr).$$
\end{lemma}

Proof: The left hand side
\begin{equation*}
\begin{split}
=& \sum_p \frac{\log^2 p}{p^{2(1+ir)}} \Bigl(\frac{1}{1 - \frac{1}{p^{1+ir}}} \Bigr)^2 = \sum_p \frac{\log^2 p}{p^{2(1+ir)}} \sum_{k=0}^{\infty} \frac{k+1}{p^{k(1+ir)}} \\
=& \sum_p \frac{\log^2 p}{p^{2(1+ir)}} \Bigl[ \sum_{k \leq \frac{\log x}{\log p}} \frac{k+1}{p^{k(1+ir)}} + O\Bigl(\int_{\frac{\log x}{\log p}}^{\infty} \frac{u}{p^u} du \Bigr) \Bigr] \\
=& \mathop{\sum_{p, k}}_{p^k \leq x} \frac{(k+1) \log^2 p}{p^{(k+2)(1+ir)}} + O\Bigl(\frac{\log x}{x}\Bigr) = \mathop{\sum_{p, l \geq 2}}_{p^{l-2} \leq x} \frac{(l-1) \log^2 p}{p^{l(1+ir)}} + O\Bigl(\frac{\log x}{x}\Bigr) \\
=& \mathop{\sum_{p, l \geq 2}}_{p^{l} \leq x} \frac{(l-1) \log^2 p}{p^{l(1+ir)}} + O\Bigl(\sum_{l = 2}^{\frac{\log 4x}{\log 2}} l \sum_{x^{1/l} < p \leq x^{1/(l-2)}} \frac{ \log^2{p}}{p^l} \Bigl) + O\Bigl(\frac{\log x}{x}\Bigr) \\
=& \mathop{\sum_{p, l \geq 2}}_{p^{l} \leq x} \frac{(l-1) \log^2 p}{p^{l(1+ir)}} + O\Bigl(\sum_{l = 2}^{\frac{\log 4x}{\log 2}} l \int_{x^{1/l}}^{\infty} \frac{\log^2 u}{u^l} du \Bigl) + O\Bigl(\frac{\log x}{x}\Bigr) \\
=& \mathop{\sum_{p, l \geq 2}}_{p^{l} \leq x} \frac{(l-1) \log^2 p}{p^{l(1+ir)}} + O\Bigl(\frac{\log^2 x \log \log x}{x^{1/2}}\Bigr) + O\Bigl(\frac{\log x}{x}\Bigr)
\end{split}
\end{equation*}
which gives the lemma.
\begin{lemma}
\label{lemma8}
Assume {\bf RH}. For real $r$, $x > 1$ and some $0 < \epsilon < 1/4$,
$$\Bigl(\frac{\zeta'}{\zeta}\Bigr)' (1+ir) + \frac{1}{r^2} = \sum_{n \leq x} \frac{\Lambda(n) \log n}{n^{1+ir}} + \frac{x^{-ir} \log x}{ir} + \frac{1 - x^{-ir}}{r^2} + O(x^{-1/4 + \epsilon}).$$
Note: Both sides are well-defined even with $r=0$.
\end{lemma}

Proof: By Perron's formula (see Titchmarsh [\ref{T}] Lemma 3.12 for example), for some large $M > 1$,
\begin{equation*}
\begin{split}
\sum_{n \leq x} \frac{\Lambda(n) \log n}{n^{1+ir}} =& \frac{1}{2\pi i} \int_{\epsilon - iM}^{\epsilon + iM} \Bigl(\frac{\zeta'}{\zeta}\Bigr)' (s+1+ir) \frac{x^s}{s} ds + O\Bigl(\frac{x^\epsilon}{M \epsilon^2} + \frac{\log^3 x}{M}\Bigr) \\
=& \mbox{Res}_{s=0} + \mbox{Res}_{s=-ir} + \frac{1}{2\pi i} \int_{-\frac{1}{2} + \epsilon - iM}^{\frac{1}{2} + \epsilon + iM} \Bigl(\frac{\zeta'}{\zeta}\Bigr)' (s+1+ir) \frac{x^s}{s} ds \\
&+ O\Bigl(\frac{x^\epsilon}{M \epsilon^2} + \frac{\log^3 x}{M}\Bigr) + O\Bigl(\frac{x^\epsilon \log^{4A} M}{M}\Bigr)
\end{split}
\end{equation*}
where $\mbox{Res}_{s=a}$ stands for the residue of the integrand at $s=a$ and the last error term comes from estimating the horizontal line integrals at height $\pm iM$ as $\frac{1}{\zeta(\sigma + it)}$, $\zeta'(\sigma + it)$ and $\zeta''(\sigma + it)$ are all smaller than $\log^A |t|$ for some $A > 0$ when $\frac{1}{2} + \epsilon \leq \sigma \leq 2$ and $|t| > 1$ under {\bf RH}. The integral on the right hand side is
$$\ll x^{-1/2 + \epsilon} M \log^{4A} M.$$
By standard calculations with Taylor's expansion, one has
$$\mbox{Res}_{s=0} = \Bigl(\frac{\zeta'}{\zeta}\Bigr)' (1+ir), \mbox{ and } \mbox{ Res}_{s=-ir} = -\frac{x^{-ir} \log x}{ir} + \frac{x^{-ir}}{r^2}.$$
Combining all these together and choosing $M = x^{1/4}$, we have the lemma.

\bigskip

Now, we are ready to resume with the proof. By Lemma \ref{lemma7} and \ref{lemma8}, the Dirichlet polynomials in $J_2$ of (\ref{3.3}) cancel out exactly and we have $J_2$
\begin{equation}
\begin{split}
=& \frac{T}{2\pi^2 \lambda} \int_{-\infty}^{\infty} \hat{k}(r \lambda) \Bigl[\frac{x^{-ir} \lambda}{ir} + \frac{1 - x^{-ir}}{r^2} + O(x^{-1/4 + \epsilon}) \Bigr] dr \\
=& \frac{T}{2\pi^2 \lambda} \int_{0}^{\infty} \hat{k}(r \lambda) \Bigl[\frac{(x^{-ir} - x^{ir}) \lambda}{ir} + \frac{2 - (x^{-ir} + x^{ir})}{r^2} \Bigr] dr + O\Bigl(\frac{T}{x^{1/4-\epsilon}}\Bigr) \\
=& \frac{T}{2\pi^2 \lambda} \int_{0}^{\infty} 2 \hat{k}(r \lambda) \Bigl[\frac{1 - \cos{(r \lambda)}}{r^2} - \frac{\sin{(r \lambda)}}{r} \lambda \Bigr] dr + O\Bigl(\frac{T}{x^{1/4-\epsilon}}\Bigr) \\
=& \frac{T}{2\pi^2} \int_{-\infty}^{\infty} \hat{k}(u) \Bigl[\frac{1 - \cos{u}}{u^2} - \frac{\sin{u}}{u} \Bigr] du + O\Bigl(\frac{T}{x^{1/4-\epsilon}}\Bigr) \\
=& \frac{T}{2\pi^2} \Bigl[ \int_{-\infty}^{\infty} \hat{k}(u) \Bigl(\frac{\sin{\frac{u}{2}}}{\frac{u}{2}}\Bigr)^2 \frac{1}{2} du - \int_{-\infty}^{\infty} \hat{k}(u) \frac{\sin u}{u} du \Bigr] + O\Bigl(\frac{T}{x^{1/4-\epsilon}}\Bigr) \\
=& \frac{T}{2\pi^2} \Bigl[ \int_{-\infty}^{\infty} k(v) 2\pi \max\{1-|2\pi v|, 0\} \frac{1}{2} dv - \int_{-\infty}^{\infty} k(v) \pi \chi_{[-\frac{1}{2}, \frac{1}{2}]} (\pi v) dv \Bigr] + O\Bigl(\frac{T}{x^{1/4-\epsilon}}\Bigr) \\
=& \frac{T}{2\pi^2} \Bigl[ 2\pi \int_{0}^{1/2\pi} k(v) (1 - 2\pi v) dv - \pi \int_{-1/2\pi}^{1/2\pi} k(v) dv \Bigr] + O\Bigl(\frac{T}{x^{1/4-\epsilon}}\Bigr) \\
=& \label{3.4} - 2T \int_{0}^{1/2\pi} k(v) v \; dv + O\Bigl(\frac{T}{x^{1/4-\epsilon}}\Bigr)
\end{split}
\end{equation}
by Parseval's identity. This cancel out with one of the terms in $\frac{1}{\lambda} \Sigma_1$. Therefore, from (\ref{3.1}) and (\ref{3.2}),
\begin{equation*}
\begin{split}
\int_{0}^{T} |S(t)|^2 dt =& R - G(T) - H(T) + O(T^{1/2} x^{1/2}) \\
=& \frac{1}{\lambda} [\Sigma_1 + \Sigma_2 + \Sigma_3] - G(T) - H(T) + O(T^{1/2} x^{1/2}).
\end{split}
\end{equation*}
Combining the results in (\ref{3.2.5}), (\ref{3.3}) and (\ref{3.4}), we have
\begin{equation}
\label{almost}
\begin{split}
& \int_{0}^{T} |S(t)|^2 dt = \frac{T}{2\pi^2} \log \log \frac{T}{2\pi} + \frac{T}{2\pi^2} \Bigl[1 + C_0 - \sum_{m=2}^{\infty} \sum_p \Bigl(\frac{1}{m} - \frac{1}{m^2}\Bigr) \frac{1}{p^m} \Bigr] \\
&- \frac{1}{\pi} \mbox{Li}\Bigl(\frac{T}{2\pi}\Bigr)
+ \frac{1}{2\pi^2 \lambda} \int_{x}^{T} \int_{-\infty}^{\infty} \hat{k}(r \lambda) \Bigl[\zeta(1-ir) \zeta(1+ir) A(ir) - \frac{1}{r^2} \Bigr] e^{-i r l} dr \; dt \\
&+ O\Bigl(\frac{T}{x^{1/4-\epsilon}}\Bigr) + O(T^{1/2} x^{1/2}) + O(x^{2 + \epsilon}) \\
=& \int_{0}^{T} \log \log \frac{t}{2\pi} + \Bigl[1 + C_0 - \sum_{m=2}^{\infty} \sum_p \Bigl(\frac{1}{m} - \frac{1}{m^2}\Bigr) \frac{1}{p^m} \Bigr] dt \\
&+ \frac{1}{2\pi^2 \lambda} \int_{x}^{T} \int_{-\infty}^{\infty} \hat{k}(r \lambda) \Bigl[\zeta(1-ir) \zeta(1+ir) A(ir) - \frac{1}{r^2} \Bigr] e^{-i r l} dr \; dt \\
&+ O\Bigl(\frac{T}{x^{1/4-\epsilon}}\Bigr) + O(T^{1/2} x^{1/2}) + O(x^{2 + \epsilon})
\end{split}
\end{equation}
\section{Completion of the proof}
It remains to deal with the integral in (\ref{almost}). Let
$$j(r) := \Bigl[\zeta(1-ir) \zeta(1+ir) A(ir) - \frac{1}{r^2}\Bigr].$$
Then $j(r) \ll \min (1, \log^B {r})$ for some $B>0$. Moreover, from the definition of $A(\eta)$, $j(z)$ can be treated as a function at least on the complex disc $|z| \leq \frac{1}{3}$. Let
$$I := \frac{1}{\lambda} \int_{-\infty}^{\infty} \hat{k}(r \lambda) j(r) e^{-irl} dr.$$
By substituting $v = \frac{rl}{2\pi}$, we have
\begin{equation}
\label{I}
\begin{split}
I =& \frac{2\pi}{\lambda l} \int_{-\infty}^{\infty} \hat{k}\Bigl(\frac{2\pi \lambda v}{l}\Bigr) j\Bigl(\frac{2 \pi v}{l}\Bigr) e^{-2\pi i v} dv \\
=& \frac{2\pi}{\lambda l} \int_{-l/8\pi}^{l/8\pi} \hat{k}\Bigl(\frac{2\pi \lambda v}{l}\Bigr) j\Bigl(\frac{2 \pi v}{l}\Bigr) e^{-2\pi i v} dv + O(e^{-c|\lambda|^{\alpha}})
\end{split}
\end{equation}
for some constant $c > 1$. Now, as $k$ and $\hat{k}$ are even, $\frac{1}{a} k(\frac{u}{a}) = \int_{-\infty}^{\infty} \hat{k}(a x) e^{-2\pi i u x} dx$. Differentiating it $m$ times, we have
\begin{equation}
\label{fourier}
\frac{1}{(-2\pi i)^m} \frac{1}{a^{m+1}} k^{(m)} \Bigl(\frac{u}{a}\Bigr) = \int_{-\infty}^{\infty} \hat{k}(a x) x^m e^{-2 \pi i u x} dx.
\end{equation}
We can do this because $\hat{k}$ has good decay. Suppose the Taylor's expansion of $j(z)$ with remainder term in $|z| \leq \frac{1}{4}$ is
\begin{equation}
\label{taylor}
a_0 + a_1 z + a_2 z^2 + ... + a_N z^N + \frac{1}{2 \pi i} \int_{C} \frac{j(s) \; ds}{s^{N+1} (s-z)} z^{N+1}
\end{equation}
where $C = \{|z| = \frac{1}{3} \}$ (see Ahlfors [\ref{A}] page 179 for example). Since $j$ is bounded on $C$, the above remainder term is $O((3|z|)^{N+1})$. Putting
(\ref{taylor}) into (\ref{I}), we have
\begin{equation*}
\begin{split}
I =& \frac{2\pi}{\lambda l} \int_{-l/8\pi}^{l/8\pi} \hat{k}\Bigl(\frac{2\pi \lambda v}{l}\Bigr) \Bigl[\sum_{n=0}^{N} a_n v^n (\frac{2\pi}{l})^n + O\Bigl(\frac{3^{N+1} v^{N+1} }{(l/2\pi)^{N+1}} \Bigr) \Bigr] e^{-2\pi i v} dv + O(e^{-c \lambda^{\alpha}}) \\
=& \frac{2\pi}{\lambda l} \sum_{n=0}^{N} a_n (\frac{2\pi}{l})^n \int_{-l/8\pi}^{l/8\pi} \hat{k}\Bigl(\frac{2\pi \lambda v}{l}\Bigr) v^n e^{-2\pi i v} dv + O\Bigl(\frac{(6\pi)^{N+1}}{l^{N+1}} \Bigr) + O(e^{-c \lambda^{\alpha}}) \\
=& \frac{2\pi}{\lambda l} \sum_{n=0}^{N} a_n (\frac{2\pi}{l})^n \int_{-\infty}^{\infty} \hat{k}\Bigl(\frac{2\pi \lambda v}{l}\Bigr) v^n e^{-2\pi i v} dv + O\Bigl(\frac{(6\pi)^{N+1}}{L^{N+1}} \Bigr) + O(e^{-d \lambda^{\alpha}})
\end{split}
\end{equation*}
for some $0 < d < 1$. Note: $a_n = \frac{1}{2\pi i} \int_{C} \frac{j(z)}{z^{n+1}} dz \ll 3^n$. By (\ref{fourier}),
\begin{equation*}
I = \frac{2\pi}{\lambda l} \sum_{n=0}^{N} a_n (\frac{2\pi}{l})^n \frac{1}{(-2\pi i)^n} \Bigl(\frac{l}{2\pi \lambda}\Bigr)^{n+1} k^{(n)}\Bigl(\frac{l}{2\pi \lambda}\Bigr) + O\Bigl(\frac{(6\pi)^{N+1}}{l^{N+1}} \Bigr) + O(e^{-d \lambda^{\alpha}}).
\end{equation*}
But recall $k(u) = -\hat{h}(u)^2 = \frac{1}{4 \pi^2 u^2}$ when $u > \frac{1}{2\pi}$. Thus,
\begin{equation}
\label{II}
\begin{split}
I =& \frac{2\pi}{\lambda l} \sum_{n=0}^{N} a_n (\frac{2\pi}{l})^n \frac{1}{(-2\pi i)^n} \Bigl(\frac{l}{2\pi \lambda}\Bigr)^{n+1} \frac{(-1)^n}{4 \pi^2} (n+1)! \Bigl(\frac{2\pi \lambda}{l}\Bigr)^{n+2} \\
&+ O\Bigl(\frac{(6\pi)^{N+1}}{l^{N+1}} \Bigr) + O(e^{-d \lambda^{\alpha}}) \\
=& \frac{1}{l^2} \sum_{n=0}^{N} a_n (\frac{-i}{l})^n (n+1)! + O\Bigl(\frac{(6\pi)^{N+1}}{l^{N+1}} \Bigr) + O(e^{-d \lambda^{\alpha}})
\end{split}
\end{equation}
which is independent of $\lambda$ and gives an asymptotic series for the lower order terms!

Apply $\Gamma(s) = \int_{0}^{\infty} e^{-u} u^{s-1} dt$ to (\ref{II}), we have, for $N \leq L^{\alpha}$,
\begin{equation*}
\begin{split}
I =& \frac{1}{l^2} \int_{0}^{\infty} \Bigl[\sum_{n=0}^{N} a_n \Bigl(\frac{-i u}{l}\Bigr)^n \Bigr] u e^{-u} du + O\Bigl(\frac{(6\pi)^{N+1}}{l^{N+1}} \Bigr) + O(e^{-d \lambda^{\alpha}}) \\
=& \frac{1}{l^2} \int_{0}^{l^{\alpha}} \Bigl[\sum_{n=0}^{N} a_n \Bigl(\frac{-i u}{l}\Bigr)^n \Bigr] u e^{-u} du + O\Bigl(\frac{(6\pi)^{N+1}}{l^{N+1}} \Bigr) + O(3^N l^{(N+1)\alpha} e^{-d \lambda^{\alpha}}) \\
=& \frac{1}{l^2} \int_{0}^{l^{\alpha}} \Bigl[ j\Bigl(\frac{-iu}{l}\Bigr) + O\Bigl(\frac{3^{N+1} u^{N+1}}{l^{N+1}}\Bigr) \Bigr] u e^{-u} du \\
&+ O\Bigl(\frac{(6\pi)^{N+1}}{l^{N+1}} \Bigr) + O(3^N l^{(N+1)\alpha} e^{-d \lambda^{\alpha}}) \\
=& \frac{1}{l^2} \int_{0}^{l^{\alpha}} j\Bigl(\frac{-iu}{l}\Bigr) u e^{-u} du + O\Bigl(\frac{3^{N+1} (N+2)!}{l^{N+1}}\Bigr) + O(3^N l^{(N+1)\alpha} e^{-d \lambda^{\alpha}})
\end{split}
\end{equation*}
Take $N = \frac{l^{\alpha}}{K \log l}$ with $K$ large enough and use $(N+2)! \leq (N+2)^{N+2}$ , we have, as $l \geq \lambda$,
\begin{equation}
\label{finish}
\begin{split}
I =& \frac{1}{l^2} \int_{0}^{l^{\alpha}} j\Bigl(\frac{-iu}{l}\Bigr) u e^{-u} du + O(e^{- d' \lambda^{\alpha}}) \\
=& \frac{1}{l^2} \int_{0}^{l^{\alpha}} u e^{-u} \Bigl[ \zeta\Bigl(1 - \frac{u}{l}\Bigr) \zeta\Bigl(1 + \frac{u}{l}\Bigr) A\Bigl(\frac{u}{l}\Bigr) + \frac{1}{(\frac{u}{l})^2} \Bigr] du + O(e^{- d' \lambda^{\alpha}}) \\
=& \int_{0}^{1} v e^{-l v} \Bigl[ \zeta(1 - v) \zeta(1 + v) A(v) + \frac{1}{v^2} \Bigr] dv + O(e^{- d' \lambda^{\alpha}})
\end{split}
\end{equation}
for some $0 < d' < 1$. Putting (\ref{finish}) into (\ref{almost}) and integrating with respect to $t$, we have Theorem \ref{theorem}. Note: $d'$ may depend on $\alpha$.
\section{Numerical evidence}
Recall $C_0 = 0.5772156649...$ and by Mathematica, one has
$$\sum_{m=2}^{\infty} \sum_{p} \Bigl(\frac{1}{m} - \frac{1}{m^2}\Bigr) \frac{1}{p^m} \approx 0.176248.$$ Note that
$$N(T) = \frac{T}{2\pi} \log{\frac{T}{2\pi e}} + \frac{7}{8} + S(T) + O\Bigl(\frac{1}{T}\Bigr)$$
where $N(T)$ denotes the number of non-trivial zeros of $\zeta(s)$ from height $0$ up to height $T$. The second moment $\int_{0}^{T} S(t)^2 dt$ is essentially the variance of $N(T)$ from $\frac{T}{2\pi} \log{\frac{T}{2\pi e}} + \frac{7}{8}$ with small error. Thus, one can use a list of the zeros of $\zeta(s)$ (for example by Odlyzko [\ref{O}]) to compute $\int_{0}^{T} S(t)^2 dt$. The asymptotic formula in (\ref{main}) can be calculated by Mathematica. For our computations, we take $x = \sqrt{T}$ and approximate the Euler product $A(r)$ by a partial product where $p$ ranges through the first $5000$ primes. Also, instead of integrating the second integral from $0$ to $1$, we integrate it from $0$ to $1 - 1/T$ with an error at most $O(\log T)$. We give the numerical evidence in the following table. Column A stands for the value from the formula in (\ref{main}), column B stands for the value from the formula in (\ref{Goldston}) and column C stands for the value from the formula in (\ref{Goldston}) subtracting $\frac{1}{\pi} \mbox{Li}(\frac{T}{2\pi})$.

\bigskip

\begin{math}
\begin{array}{l|l|l|l|l}
T & \int_{0}^{T} S(t)^2 dt & A & B & C \\
\hline
9998.85040 & 1653.145 & 1651.05 & 1721.61 & 1638.76 \\
19999.27562 & 3411.009 & 3407.72 & 3534.54 & 3386.35 \\
29999.71003 & 5200.768 & 5196.71 & 5376.49 & 5167.12 \\
39999.49733 & 7009.117 & 7005.47 & 7236.31 & 6968.18 \\
49999.57275 & 8831.813 & 8828.58 & 9109.15 & 8783.93 \\
59998.88155 & 10668.969 & 10662.7 & 10991.9 & 10610.9 \\
69999.61050 & 12509.875 & 12506.1 & 12883.3 & 12447.4
\end{array}
\end{math}

\bigskip

The values of $T$ are the largest imaginary parts of the zeros of $\zeta(s)$ just below $10000$, $20000$, ..., $70000$ respectively. By comparing the second and third column, the asymptotic formula in (\ref{main}) gives very good approximation to the second moment. By comparing the second, fourth and last column, we see that it is reasonable to have $-\frac{1}{\pi} \mbox{Li}(\frac{T}{2\pi})$ in the asymptotic formula.

Tsz Ho Chan\\
American Institute of Mathematics\\
360 Portage Avenue\\
Palo Alto, CA 94306\\
USA\\
thchan@aimath.org

\end{document}